\documentclass{amsart}

\usepackage{amssymb}

\newtheorem{thm}{Theorem}

\newtheorem{lem}[thm]{Lemma}
\newtheorem{cor}[thm]{Corollary}

\newtheorem{prop}[thm]{Proposition}

\theoremstyle{definition}
\newtheorem{defn}[thm]{Definition}

\newtheorem{say}[thm]{}
\newtheorem{exmp}[thm]{Example}

\newtheorem{const}[thm]{Construction}   

\newtheorem{rem}[thm]{Remark}          
\newtheorem{ack}{Acknowledgments}

\newtheorem{defn-thm}[thm]{Definition--Theorem}  
\newtheorem{defn-lem}[thm]{Definition--Lemma}  

\theoremstyle{remark}

\renewcommand{\c}[0]{{\mathbb C}}  

\renewcommand{\o}[0]{{\mathcal O}} 
\newcommand{\z}[0]{{\mathbb Z}}

\renewcommand{\r}[0]{{\mathbb R}}

\newcommand{\p}[0]{{\mathbb P}}

\newcommand{\q}[0]{{\mathbb Q}}

\newcommand{\qtq}[1]{\quad\mbox{#1}\quad}

\newcommand{\rank}[0]{\operatorname{rank}}

\newcommand{\supp}[0]{\operatorname{Supp}}

\newcommand{\im}[0]{\operatorname{im}}

\newcommand{\coker}[0]{\operatorname{coker}}

\newcommand{\onto}[0]{\twoheadrightarrow}

\newcommand{\lcm}[0]{\operatorname{lcm}}

\newcommand{\stab}[0]{\operatorname{Stab}}

\def\into{\DOTSB\lhook\joinrel\to}

\begin{document}
\bibliographystyle{amsalpha}

\title{Circle actions on simply connected 5--manifolds}
\author{J\'anos Koll\'ar}

\today
\maketitle
\begin{abstract} 
The aim of this paper is to study compact 5--manifolds which
admit
 fixed point free circle actions.
The first result implies that 
 the torsion in the second homology 
and the second Stiefel--Whitney class
have to satisfy strong restrictions.
We then show that for simply connected 5--manifolds
these restrictions
are necessary and sufficient.
\end{abstract}

\tableofcontents

It is easy to see  that a
simply connected compact 5--manifold $L$ admits a free circle action
iff $H_2(L,\z)$ is torsion free and
the classification of free circle actions
up to diffeomorphism is equivalent
to the classification of 
simply connected compact 4--manifolds plus the action of their
diffeomorphism group on the 
second cohomology (cf.\ \cite[Prop.10]{geig}).

Motivated by some questions that arose in connection
with the study of complex analytic Seifert $\c^*$-bundles
\cite{ko-s2s3, ko-es5}, 
 this paper  investigates  compact 5-manifolds that admit
circle actions where the stabilizer of every point is finite,
that is, fixed point free circle actions.
We show that in this case  $H_2(L,\z)$  can contain torsion,
but the torsion and the second Stiefel--Whitney class
have to satisfy strong restrictions.
We then show that for simply connected manifolds
these restrictions
are necessary and sufficient for the existence of a
fixed point free circle action.

\begin{defn}\label{i(L).defn}
Let $M$ be any manifold. Write its  second homology as
 a direct sum
of cyclic groups of prime power order
$$
H_2(M,\z)=\z^k+\sum_{p,i} \bigl(\z/p^i\bigr)^{c(p^i)}
\qtq{for some $k=k(M), c(p^i)=c(p^i,M)$.}
\eqno{(\ref{i(L).defn}.1)}
$$
The numbers $k, c(p^i)$ are  determined by
$H_2(M,\z)$ but the subgroups
$(\z/p^i)^{c(p^i)}\subset H_2(M,\z)$ are  usually not unique.
One can choose the decomposition
(\ref{i(L).defn}.1) such that 
the second Stiefel--Whitney class map
$$
w_2:H_2(M,\z)\to \z/2
$$ is zero 
on all but one summand  $\z/2^n$. This value $n$ is unique
and it is   denoted by $i(M)$ \cite{barden}.
This invariant
 can take up any value $n$ for which $c(2^n)\neq 0$,
besides $0$ and $\infty$. 
Alternatively, 
$i(M)$ is the
smallest $n$ such that there is an $\alpha\in H_2(M,\z)$ such that
$w_2(\alpha)\neq 0$ and $\alpha$ has order $2^n$.
\end{defn}

The existence of a 
fixed point free differentiable circle action 
puts strong restrictions on $H_2$ and on $w_2$.

\begin{thm}\label{main.thm}
Let  $L$ be a   compact  5--manifold with $H_1(L,\z)=0$
which admits a fixed point free differentiable circle action.
Then:
\begin{enumerate}
\item  For every prime $p$, we have
at most $k+1$ nonzero $c(p^i)$ in 
(\ref{i(L).defn}.1). That is,  
$\#\{i:c(p^i)>0\}\leq k+1$.
\item One can arrange that $w_2:H_2(L,\z)\to \z/2$
is the zero map on all but the $\z^k+(\z/2)^{c(2)}$ summands
in (\ref{main.thm}.3). That is, $i(L)\in \{0,1,\infty\}$.
\item If $i(L)=\infty$ then $\#\{i:c(2^i)>0\}\leq k$.
\end{enumerate}
\end{thm}

These conditions are sufficient
for simply connected manifolds:

\begin{thm}\label{main.thm.2}
Let  $L$ be a   compact, simply connected  5--manifold.
Then $L$ admits a fixed point free differentiable circle action
if and only if $w_2:H_2(L,\z)\to \z/2$ satisfies
the conditions  (\ref{main.thm}.1--3).
\end{thm}

The conditions are especially transparent for
homology spheres.

\begin{exmp} Let $c(p^i)$ be any sequence of even natural numbers,
only finitely many nonzero. By \cite{smale}, there is a unique
simply connected, spin, compact  5--manifold $L$ such that
$H_2(L,\z)\cong \sum_{p,i}\bigl(\z/p^i\bigr)^{c(p^i)}$.

By (\ref{main.thm}), this $L$ admits
 a fixed point free differentiable circle action
iff for every $p$, at most one of the
$c(p), c(p^2), c(p^3),\dots$ is nonzero.
\end{exmp}

It should be noted that the proof does {\em not} give a classification of
all  fixed point free  
$S^1$-actions on any 5--manifold. In fact, we exhibit
infinitely many topologically distinct  fixed point free  $S^1$-actions
on every $L$ as in (\ref{main.thm.2}). In principle the 
classification of
all $S^1$-actions on 5--manifolds is reduced to a
question on 4--dimensional orbifolds, but the 4--dimensional
question is rather complicated.

\begin{say}\label{or-wa.say}
The classification of fixed point free circle actions on 3--manifolds
was considered by Seifert \cite{seif}. If $M$ is a  3--manifold
with a fixed point free circle action then the quotient space
$F:=M/S^1$ is a surface (without boundary in the orientable case).
The classification of these {\it Seifert fibered} 
3--manifolds $f:M\to F$ is thus equivalent to the
 classification of fixed point free circle actions.
It should be noted that already in this classical case,
it is conceptually better to view the base surface $F$
not as a 2--manifold but as a 2-dimensional {\it orbifold},
see \cite{scott} for a detailed survey from this point of view.

The classification of circle actions on
  4--manifolds is treated in
\cite{fint1, fint2}.  Here the quotient is a 3--manifold
(with boundary corresponding to the fixed points) endowed
with additional data involving links and  certain weights.

A generalization of Seifert bundles to higher dimensions
was considered in \cite{or-wa}. In essence, this paper considers the case
 when  $L$ is a real hypersurface
in a complex manifold $Y$ with a $\c^*$-action,
$Y/\c^*$ is a complex manifold
and $L$ is invariant under the induced $S^1$-action. 
The computations of \cite{or-wa} are topological in nature, and use only
that  $Y/\c^*=L/S^1$ is a real manifold, and
  the fixed point set of every element of
$S^1$ is oriented.
These assumptions do not hold in general, 
and we   follow a somewhat different approach.

Foundational  questions concerning circle actions are also considered in
\cite{HaSa91}.
\end{say}

The proofs of (\ref{main.thm}) and (\ref{main.thm.2}) follow the path
of \cite{seif, or-wa}. We 
start with any manifold $L$ with a  fixed point free  circle action, and
consider the
quotient space $X:=L/S^1$. $X$ is usually not a manifold,
only an orbifold, but we consider it with a much richer orbifold structure
$(X,\Delta)$ where $\Delta=\sum(1-\frac1{m_i})D_i$  is a formal sum
of codimension 2 closed subspaces $D_i\subset X$.
The main technical aspect of the proof is to understand how to
relate invariants of 
$L$ and invariants of the  orbifold $(X,\Delta)$.

In order to prove 
(\ref{main.thm}),
we try to compute the Leray spectral sequence
$$
H^i(X, R^jf_*\z_L)\Rightarrow H^{i+j}(L,\z).
$$
This is similar to the Gysin sequence
used in \cite{or-wa}, but the Leray spectral sequence
is better suited to the current situation.
We end up computing $H_2(L,\z)$ in terms of
the  orbifold $(X,\Delta)$, but the formula involves
$H^{\dim X-3}(X,\z)$ which I can not control in general.
If $\dim L=5$ then this is $H^1(X,\z)$ and it 
vanishes if $H_1(L,\z)=0$.

To see that the restrictions of (\ref{main.thm}) are
sufficient,
 we provide examples of Seifert bundles
$L\to X$ with $X=(k+1)\#\c\p^2$, a connected sum of $k+1$ copies of
$\c\p^2$, and the $D_i\subset X$ are smooth  surfaces
intersecting transversally. 
It is somewhat lucky that these special cases cover all
possibilities.
For these examples 
we compute $\pi_1(L)$, $w_2(L)$ and 
$H_2(L,\z)$. Everything is easier since we do not have to
worry about  orbifold points of $X$.
We then conclude the proof by using
 the structure theorem of 
simply connected compact 5--manifolds due to
Smale and Barden.

\begin{thm} \cite{smale, barden}\label{smale.thm}
Let $L$ be a  simply connected compact 5--manifold.
Then $L$ is uniquely determined by 
$H_2(L,\z)$ and the second
Stiefel--Whitney class map
$w_2: H_2(L,\z)\to \z/2$.

Furthermore, there is such a 5--manifold 
iff there is an integer  $k\geq 0$ and a finite Abelian group $A$
such that either
\begin{enumerate}
\item $H_2(L,\z)\cong \z^k+A+A$ and
$w_2: H_2(L,\z)\to \z/2$ is arbitrary, or 
\item $H_2(L,\z)\cong \z^k+A+A+\z/2$ and
$w_2$ is projection on the $\z/2$-summand.
\end{enumerate}
\end{thm}

My original interest in this topic came from
complex geometry. A method of Kobayashi \cite{kob}, generalized in 
\cite{bg00, bgn03c},
allows one to construct positive Ricci curvature 
Einstein metrics on $L$ from a
 positive Ricci curvature  orbifold K\"ahler--Einstein  metric on $(X,\Delta)$ 
if the base orbifold $(X,\Delta)$
has a  complex structure.
The existence of a positive Ricci curvature orbifold 
  K\"ahler--Einstein metric on $(X,\Delta)$ 
imposes strong restrictions. These were explored in
\cite{ko-s2s3, ko-es5}. It seemed to me, however, that behind the
conditions coming from complex geometry, there were some
weaker but non obvious topological restrictions as well.

In fact, the $k=0$ case of
(\ref{main.thm}.1) first appeared in \cite[Cor.81]{ko-es5}
as a restriction on  
Seifert bundles over algebraic orbifolds.
Now we see that this restriction is imposed not by the
presence of an algebraic structure but by the
topological circle action. On the other hand, in
 (\ref{sasak.rests}) we exhibit  additional,
this time non--topological, restrictions on
Seifert bundles over algebraic orbifolds.

\section{Local classification of $S^1$-actions}

\begin{defn}\label{action.defn}
Let $M$ be a differentiable manifold with a differentiable circle action
$\sigma:S^1\times M\to M$. 
I usually think of $S^1$ as a subgroup of $\c^*$.
Pick a point $p\in M$ which is not a fixed point
and let $O(p)\subset M$ be the orbit of $p$. The stabilizer of $p$,
denoted by 
$\stab_p\subset S^1$, is cyclic of order $m=m(p)$ and we can choose a 
canonical generator  $e^{2\pi i/m}\in \stab_p$.

Let $p\in H_p$ be a codimension 1 submanifold, transversal to
$O(p)$, invariant under $\stab_p$. 
Let $T_p$ denote the tangent space of $H_p$ at $p$
with its induced faithful $\stab_p$-action. 
This action of $\stab_p$ on $T_p$ is  the
{\it stabilizer} or {\it slice} representation.
A neighborhood of
$O(p)$ in $M$ is $S^1$-equivariantly diffeomorphic to
$$
S^1\times T_p/\stab_p,\qtq{where}
\eqno{(\ref{action.defn}.1)}
$$
\begin{enumerate}
\item[i)]  the $S^1$-action is the natural
$S^1$-action on itself, and
\item[ii)]  $\stab_p$ acts on $S^1$ by multiplication 
as a subgroup and on $T_p$
by the {\em inverse} of the stabilizer representation.
\end{enumerate}
Thus the local structure of $\sigma:S^1\times M\to M$
near any orbit is determined by the stabilizer
representation of $\stab_p$ on $T_p$.

If the action 
$\sigma:S^1\times M\to M$ has no fixed points,
set $X:=M/S^1$ and let $f:M\to X$ denote the quotient map.
We call $f:M\to X$ the  {\it Seifert bundle} 
associated to the circle action.
Later we  modify this definition slightly
and view $X$ not as a topological space but as a 
differentiable  orbifold.

Every fiber of $f$ is an $S^1$-orbit, thus a circle.
For $x\in X$, the stabilizer $\stab_p$ is independent of $p\in f^{-1}(x)$.
Its order is called the {\it multiplicity} of the
fiber $f^{-1}(x)$ and it is denoted by $m(x)$ or
$m(x,M)$ if there is some doubt as to which $M$ we work with.

Given $x\in X$ and $p\in f^{-1}(x)$,
an open neighborhood of $x$ is also realized as
$T_p/\stab p$. This gives $X$ the structure of a
{\it cyclic orbifold}. That is, it is patched together from
orbifold charts of the form $\r^n/(\mbox{cyclic group})$
by a linear action. 

We give a detailed local description of the orbifold
structure on $M/S^1$ later.
\end{defn}

\begin{defn} We say that $\sigma:S^1\times M\to M$
has  {\it orientable stabilizer representations} if 
 the representation of $\stab_p$ on $T_p$
is orientation preserving for every $p$. 

If $M$ itself is orientable, 
then every element of a connected group acting on $M$
preserves orientation, hence the
stabilizer representations are all
orientable.

This condition is
 also satisfied in many other cases when $M$ is not orientable.
The nonorientable  stabilizer representation case is anomalous already
in dimension 3 \cite{scott}.
\end{defn}

\begin{say}[Real representations of cyclic groups]
\label{real.orient.reps}
The  orientation preserving irreducible real representations
of a cyclic group $\z/m$ are the trivial representation $R_{m,0}$ 
and the 2--dimensional representations
$$
R_{m,j}: \z/m\ni 1\mapsto
\left(
\begin{array}{cc}
\cos \tfrac{2\pi j}{m} & \sin \tfrac{2\pi j}{m}\\[1ex]
-\sin \tfrac{2\pi j}{m} & \cos \tfrac{2\pi j}{m}
\end{array}
\right),\quad j=1,\dots, m-1.
$$
If $m$ is even then $R_{m,m/2}$ decomposes as the sum
of two orientation reversing representations.

If $V$ is a 2--dimensional  faithful
representation of $\z/m$, then either
$V$ is orientation preserving or
$m=2$.

$R_{m,j}$ is orientation reversing isomorphic to
$R_{m,m-j}$. A faithful irreducible real representation
 is orientation reversing isomorphic to
itself only for $m=2$.

Any orientation preserving real representation
of $\z/m$ can  be written as
the direct sum of orientation preserving  irreducible  real representations.
(This is not quite unique as $R_{m,j}+R_{m,j}$ and $R_{m,m-j}+R_{m,m-j}$
are orientation preserving isomorphic.)
Thus every orientation preserving real representation
on $\r^{2n}$ can be obtained from
a complex  representation on $\c^n$ by forgetting the complex structure.
This correspondence is, however, not entirely natural, as
we need to specify an orientation on each
irreducible subrepresentation.
\end{say}

From (\ref{action.defn}.1) we obtain the following:

\begin{lem} \label{cod.rep.lem}
Let  $M$ be a differentiable manifold with a circle action
$\sigma:S^1\times M\to M$ with orientable stabilizer representations.
Given integers $m$ and $c_1,\dots, c_{m-1}$, the set of all points
$M^0(m,c_1,\dots, c_{m-1})\subset M$ where
$\stab_p=\z/m$ and the representation of $\stab_p$ on $T_p$ is
isomorphic to
$$
\sum c_jR_{m,j}+(\mbox{trivial representation})
$$
is a smooth 
submanifold of codimension $2\sum c_j$ (or empty).\qed
\end{lem}

\begin{say}[Codimension 2 fixed points]  \label{cod.2.say}
By (\ref{cod.rep.lem}), the codimension 2 fixed points
correspond to stabilizer representations
$$
R_{m,j}+(\mbox{trivial representation}).
$$
Let us denote the corresponding subset of $M$ by $M^0(m,j)$.
 Note that for now we have a some
nonuniqueness since we can not distinguish $M^0(m,j)$ from $M^0(m,m-j)$.
This will be rectified later by fixing some orientations.
(This notation gives two possible meanings to $M^0(2,1)$,
but we end up with the same submanifold.)

Fix a point $p\in M^0(m,j)$. Depending on the  
orientation of $T_p$,
the stabilizer representation is $R_{m,j}$ or $R_{m,m-j}$.
If $m\geq 3$, then these are not orientation preserving isomorphic,
so fixing say $R_{m,j}$ gives a well defined  orientation of $T_p$.

As we move the point $p$ in a connected component of
$M^0(m,j)$, we get an orientation 
of $T_p$ for every $p$. Thus we obtain:

\begin{lem}\label{norm.bund.orient} If $m\geq 3$, the normal bundle
of $M^0(m,j)$ in $M$ is orientable.\qed
\end{lem}

{\it Note.} If $\dim M=3$ then each connected component of 
$M^0(m,j)$ is a single $S^1$-orbit,
hence naturally oriented. In the  cases considered by \cite{or-wa},
the normal bundle to $M^0(m,j)$ has a complex structure
hence  a natural orientation.
It seems to me that in general there is no natural choice,
and these orientations have to be chosen by hand.
\end{say}

If $T_p$ has dimension $n$ and the 
corresponding  stabilizer representation is
$R_{m,j}+(\mbox{trivial representation})$, then
we can write 
$$
T_p\cong \c_z+\r^{n-2},
$$
and the representation  is given by
multiplication by $\epsilon^j$ on $\c_z$ where
$\epsilon=e^{2\pi i/m}$ and $(j,m)=1$. Thus we can write
$$
T_p/\stab_p\cong \c_x+\r^{n-2}\qtq{where} x=z^m.
$$
In particular, at these points the quotient $M/S^1$ is a manifold.

\begin{defn} 
Let  $M$ be a differentiable manifold with a circle action
$\sigma:S^1\times M\to M$ with orientable stabilizer representations.
Let $M^s\subset M$ be the closed  set of points $p\in M$ where the
invariant subspace of the stabilizer representation
has codimension at least 4.

Set $X=M/S^1$  with quotient map $f:M\to X$ and $X^s=M^s/S^1$.
As noted above, $X^0:=X\setminus X^s$ is a manifold
and 
$$
f(\cup_{mj} M^0(m,j))\subset X^0
$$
is a closed submanifold of codimension 2.
Let its connected components be $D_i^0$.
We see in (\ref{cod.4.prop}.4) that their closures
$D_i\subset X$ are suborbifolds.

Each $D_i^0$ lies in the image of a unique
$M^0(m,j)$. This assigns a natural number
$m=m_i$ to $D_i$. We introduce the formal notation
$$
\bigl(X,\sum_i (1-\tfrac1{m_i})D_i\bigr)
$$
to denote the  {\it base orbifold} of
$f:M\to X$. 
We call
$$
f:M\to \bigl(X,\sum_i (1-\tfrac1{m_i})D_i\bigr)
$$
the  {\it  Seifert bundle} associated to the circle action
$\sigma:S^1\times M\to M$. Sometimes we use the shorthand
$\Delta=\sum_i (1-\tfrac1{m_i})D_i$.

(The  choice of the  coefficients $1-\tfrac1{m_i}$
 comes from complex geometry where the orbifold
first Chern class is given by the formula
$c_1(X)-\sum_i (1-\tfrac1{m_i})[D_i]$.)

We see in (\ref{cod.4.say}) that the data $(X,D_i,m_i)$
 determine the orbifold structure of $X$.

From (\ref{norm.bund.orient}) we conclude that the
normal bundle of $D_i^0\subset X^0$ is orientable
if $m_i\geq 3$. 

If $\bigl(X,\sum_i (1-\tfrac1{m_i})D_i\bigr)$
is an orbifold with  oriented normal bundles $N_{D_i,X}$ for $m_i\geq 3$
then this gives an orientation to the normal bundle
of $M^0(m,j)$ for $m\geq 3$ and this distinguishes
$M^0(m,j)$ from $M^0(m,m-j)$. 
Thus each $D_i^0$ with $m_i\geq 3$ is in the image of a unique
$M^0(m_i,j_i)$. Since the stabilizer representation is faithful,
$j_i$ is relatively prime to $m_i$, so $b_ij_i\equiv 1\mod m_i$
has a unique solution $1\leq b_i<m_i$.
If $m_i=2$ then $j_i=1$ so we can take $b_i=1$.

The pair $(m_i,b_i)$ is called the {\it orbit invariant} along $D_i$. 
(It is  denoted by $(\alpha_i,\beta_i)$ in \cite{or-wa}.)
Again I emphasize that while $m_i$ and the unordered pair
$\{b_i,m_i-b_i\}$ are determined by
$f:M\to X$, one needs an orientation of the normal bundle
of $D^0_i$ to determine $b_i$ itself.

We see in (\ref{cod.4.prop}.6) that the data
$$
\bigl(X, \textstyle{\sum_i}\tfrac{b_i}{m_i}D_i\bigr)
$$
determine $M$ locally on $X$ if $X$ is smooth. 
If $X$ is not a manifold, one needs further local invariants
at the singular points of $X$.

On the other hand,  $\bigl(X, \textstyle{\sum_i}\tfrac{b_i}{m_i}D_i\bigr)$
 does not determine
$M$ globally, but the different choices are obtained by
``twisting'' with $H^2(X,\z)$  
(\ref{seif.obstr}).
\end{defn}

\begin{say}[Codimension 4 fixed points]  \label{cod.4.say}

As we noted in (\ref{real.orient.reps}),
 an $n$-dimensional orientable real representation of
$\z/m$ can be written as
$\c^k+\r^{n-2k}$ where the action is trivial on $\r^{n-2k}$
and on $\c^k$ it is given by
$$
(z_1,\dots,z_k)\mapsto (\epsilon^{j_1}z_1,\dots,\epsilon^{j_k}z_k)
\qtq{where $\epsilon=e^{2\pi i/m}$.}
$$
The summand $\r^{n-2k}$ does not give any interesting contribution, and we
concentrate on the $\c^k$ part. As a shorthand, we denote
 the corresponding quotient by
$$
\c^k/\tfrac1{m}(j_1,\dots,j_k).\eqno{(\ref{cod.4.say}.1)}
$$
We write $\c^k_{\mathbf z}$ or $\c^k_{z_1,\dots,z_k}$  to indicate
the name of the coordinates.

The corresponding Seifert bundle is given locally by
$$
\bigl(S^1_v\times
\c^k_{z_1,\dots,z_k}\bigr)/\tfrac1{m}
(1,-j_1,\dots,-j_k)+\r^{n-2k} \to
 \c^k_{z_1,\dots,z_k}/\tfrac1{m}
(j_1,\dots,j_k)+\r^{n-2k},
$$
where the sign change is coming from (\ref{action.defn}.1.ii).
We usually drop the uninteresting $\r^{n-2k}$ in the sequel.

 It is also useful to
extend $S^1_v$ to $\c^1_v$ and consider the quotient
$$
\bigl(\c_v\times
\c^k_{z_1,\dots,z_k}\bigr)/\tfrac1{m}
(1,-j_1,\dots,-j_k) \to
 \c^k_{z_1,\dots,z_k}/\tfrac1{m}
(j_1,\dots,j_k),
$$
which is the corresponding Seifert $\c^*$-bundle
extended by the zero section.

Given $j_1,\dots,j_k$ and $m$, set
 $$
c_i:=\gcd(j_1,\dots,\widehat{j_i},\dots, j_k,m),\quad
d_i:=j_ic_i/C\qtq{and}C:=\prod c_i.
\eqno{(\ref{cod.4.say}.2)}
$$
Note that the $c_i$ are pairwise relatively prime and $C/c_i$ divides $j_i$.
Observe that $\z/c_i\subset \z/m$ acts trivially on all but the $i$th
coordinate of $\c^k$,
 so it is a quasi reflection. 

In particular, $(z_i=0)$ is the closure of the unique
connected component of $M^0(c_i,j_i)$ intersecting this chart.

The conditions 
$$
r\equiv j_i\mod c_i\qtq{and} c_{i'}|r\ \mbox{for}\ i'\neq i
$$
imply $r\equiv j_i\mod C$. Thus the
codimension 2 orbit invariants  $(c_i, b_i)$ (with $b_ij_i\equiv 1\mod c_i$)
determine  $(m, j_1,\dots,j_k)$, and hence the
local structure of the Seifert bundle, 
if $m=\prod c_i$.

The quotient
of $\c^k_{\mathbf z}$ by $\z/C\cong \sum \z/c_i$
is again an affine space $\c^k_{\mathbf x}$ with $x_i=z_i^{c_i}$.
Thus, as a topological space
$$
\c^k_{\mathbf z}/\tfrac1{m}(j_1,\dots,j_k)\cong 
\c^k_{\mathbf x}/\tfrac1{m/C}(d_1,\dots,d_k).
$$
The fixed point set of every nonidentity element of $\z/(m/C)$ has 
complex codimension $\geq 2$, thus
$\c^k_{\mathbf x}/(\z/(m/C))(d_1,\dots,d_k)$ is a manifold
only if $m/C=1$.

\end{say}

We can summarize these results as follows.

\begin{prop}\label{cod.4.prop}
Let  $M$ be a differentiable manifold with a circle action
$\sigma:S^1\times M\to M$ with orientable stabilizer representations.
Let $p\in M^s$ be a point with stabilizer representation
$$
R_{m,j_1}+\cdots+R_{m,j_k}+\r^{n-2k}.
$$
Let $c_i$ be defined as in  (\ref{cod.4.say}.2). Then
\begin{enumerate}
\item An open dense subset of $(z_i=0)$ is contained in
$M^0(c_i,j_i)$. These are the only $M^0(m,j)$ whose closure contains $p$.
\item The $c_i$ are pairwise relatively prime
and $\prod c_i$ divides $m$.
\item $M/S^1$ is a manifold at the image of $p$ iff
$m=\prod c_i$. In this case
$$
\c^k_{\mathbf z}/\tfrac1{m}(j_1,\dots,j_k)\cong 
\c^k_{\mathbf x}\qtq{with $x_i=z_i^{c_i}$.}
$$
\end{enumerate}
Translating these into global terms we get the following:
\begin{enumerate}\setcounter{enumi}{3}
\item The closures $M(m,j)$ of $M^0(m,j)\subset M$ are smooth and 
intersect each other transversally.
\item If $M(m_1,j_1)\cap M(m_2,j_2)\neq \emptyset$
then $(m_1,m_2)=1$.
\item $M/S^1$ is a manifold iff for every $p\in M$
$$
|\stab_p| =\prod_{M(m,j)\ni p} m.
$$
In this case the pairs $\{(m,j): M(m,j)\ni p\}$
determine the $S^1$-action in a neighborhood of
the orbit $O(p)$. \qed
\end{enumerate}
\end{prop}

Some special properties of the 5--dimensional case
are worth emphasizing:

\begin{say}\label{dim5.nice.say} 
Let  $L$ be a differentiable 5--manifold with a circle action
$\sigma:S^1\times L\to L$ with orientable stabilizer representations
and $f:L\to (X, \sum(1-\frac1{m_i})D_i)$ the corresponding
Seifert bundle. Then
\begin{enumerate}
\item $X^s\subset X$ is finite and $X\setminus X^s$ is a manifold,
\item each $D_i\subset X$ is a 2--manifold and at most 2 of them pass
through any point of $X$.
\end{enumerate}
\end{say}

\section{The cohomology groups  of Seifert bundles}

In working with group actions on manifolds and taking various
quotients, one frequently runs into orbifolds.
It is therefore  convenient to define the notion of
Seifert bundles in a rather general setting.
In order to avoid pointless complications, let us
assume from now on that every topological space
is a CW complex.

\begin{defn}  \label{top.seif.Gbund.defn}
A   {\it generalized Seifert bundle} over 
$X$ is a topological space
$M$ together with an 
 $S^1$-action and a
continuous  map $f:M\to X$
such that  
$X$ has an open covering $X=\cup_iU_i$
such that for every $i$ 
  the preimage $f:f^{-1}(U_i)\to U_i$ is
fiber preserving $S^1$-equivariantly homeomorphic to
a ``standard local  generalized Seifert bundle''
$$
f_i:(S^1\times V_i)/(\z/m_i)\to U_i.
$$
Here $V_i$ is a topological  space with a $\z/m_i$-action
  such that $V_i/(\z/m_i)\cong U_i$
and
 the $\z/m_i$-action on $S^1\times V_i$ is the diagonal action
given 
on $S^1$ by a homomorphism
$\phi_i:\z/m_i\to S^1$ composed with  the action of $S^1$ on itself.
 The action of $S^1$ on itself
gives the $S^1$ action on $Y$.

In order to avoid 
nontrivial orbifold structures on $M$, we always assume
that the $\z/m_i$-action on $S^1\times V_i$ is fixed point free outside
a codimension 2 set.

For $x\in X$ let $U_i\ni x$ be an open subset as above.
Let $v\in V_i$ be a preimage of $x$. 
If $\stab_v\subsetneq \z/m_i$ is a proper subgroup, then
there  are
open subsets $x\in U_x\subset U_i$ 
and $v\in V_x\subset V_i$ such that 
$U_x=V_x/\stab_v$ and we can also describe our
 generalized Seifert bundle locally as
$$
f^{-1}(U_x)\cong (S^1\times V_x)/\stab_v.
$$
The order of the group $\stab_v$ depends only on $x$,
and it is called the {\it multiplicity} of the fiber of
$f:M\to X$ over $x$. It is denoted by $m(x)$ or $m(x,M)$.
\end{defn}

\begin{say}[Maps between  generalized Seifert bundles]\label{seif.maps.say}

Let $f:M\to X$ be a  generalized Seifert bundle.
The $S^1$-equivariant homeomorphisms $h:M\to M$ such that $f\circ h=f$
are  the multiplications
$p\mapsto \phi(f(p))\cdot p$ where $\phi:X\to S^1$ is any continuous function,
cf.\ \cite[3.1]{HaSa91}.

It is more interesting to look at higher degree maps 
$h:M_1\to M_2$ between  generalized Seifert bundles.

Let $f:M\to X$ be a  generalized Seifert bundle and $\z/m\subset S^1$ a finite subgroup.
Then the $S^1$ action descends to an $S^1/(\z/m)$-action on 
$M/(\z/m)$
and  $f/(\z/m):M/(\z/m)\to X$ is another  generalized Seifert bundle.

Even if $M$ is  manifold, in general $M/(\z/m)$ is only an orbifold.

The case when $m=m(X)$ is the least common multiple of the
multiplicities $m(x)$ is especially useful. We denote this quotient by
$f/\mu:M/\mu\to X$.  Since every stabilizer $\stab_p: p\in M$
is contained in $\z/m(X)$, we conclude that 
$f/\mu:M/\mu\to X$ is a locally trivial $S^1$-bundle.
Locally trivial $S^1$-bundles are classified by their Chern class
$c_1((M/\mu)/ X)\in H^2(X,\z)$ and we define the
{\it Chern class} of the  generalized Seifert bundle $f:M\to X$ as
$$
c_1(M)=c_1(M/X):=\tfrac1{m(X)}c_1((M/\mu)/ X)\in H^2(X,\q).
$$
Usually it is not an integral cohomology class.
\end{say}

Our aim 
is to obtain information about 
 the integral cohomology groups 
of a  generalized Seifert bundle
$f:M\to (X,\Delta)$ in terms of
$(X,\Delta)$ and the Chern class  of $M/X$.

The cohomology groups $H^i(M,\z)$ are computed by a 
Leray  spectral sequence whose $E_2$ term is
$$
E^{i,j}_2=H^i(X, R^jf_*\z_{M})\Rightarrow H^{i+j}(M,\z).
$$
Every fiber  of $f$ is $S^1$, so 
$R^jf_*\z_M=0$ for $j\geq 2$ and
the only
interesting
higher direct image is $R^1f_*\z_{M}$.
Our first task is to compute this sheaf and its
cohomology groups. 
Next we consider the edge homomorphisms in the spectral
sequence
$$
\delta_i: H^i(X, R^1f_*\z_{M})\to H^{i+2}(X,\z),
$$
and identify them, at least modulo torsion,
 with cup product with  the Chern class  $c_1(M/X)$.

In some cases of interest, these data completely determine the
cohomology groups, and even the topology,  of $M$.
Some of these instances are discussed in
\cite{or-wa, ko-s2s3, ko-es5}.

\begin{prop}\label{R^1.prop} Let $f:M\to X$ be a  generalized Seifert bundle.
\begin{enumerate}
\item There is a natural isomorphism
$\tau_M:R^1f_*\q_M\cong \q_X$.
\item There is a natural injection
$\tau_M:R^1f_*\z_M\into \z_X$ 
which is an isomorphism over points where $m(x)=1$.
\item If $U\subset X$ is connected then
$$
\tau_M(H^0(U,R^1f_*\z_M))= m(U)\cdot  H^0(U,\z)\cong m(U)\cdot \z,
$$
where $m(U)$ is the $\lcm$ of the multiplicities of
all fibers over $U$.
\end{enumerate}
\end{prop}

Proof. Pick $x\in X$ and a small contractible neighborhood
$x\in V\subset X$. Then $f^{-1}(V)$ retracts to
$S^1\subset f^{-1}(x)$ and (together with the orientation of $S^1$)
this gives  a distinguished generator
$\rho\in H^1(f^{-1}(V), \z)$. This in turn determines
a  cohomology class $\frac1{m(x)}\rho\in H^1(f^{-1}(V), \q)$.
These normalized cohomology classes  are compatible with each other
and give a global section of $R^1f_*\q_M$. Thus
$R^1f_*\q_M=\q_X$ and we also obtain the
injection $\tau:R^1f_*\z_M\into \z_X$ as in (2). 

If $U\subset X$ is connected, a section $b\in \z\cong H^0(U,\z_U)$ is in
$\tau(R^1f_*\z_M)$ iff $m(x)$ divides $b$ for every $x\in U$.
This is exactly (3).\qed

\begin{say}\label{chernclass=edgemap.say}
Given $f:M\to X$, consider
 the quotient map $\pi:M\to M/\mu=M/(\z/m(X))$
 defined in (\ref{seif.maps.say}).
Apply  (\ref{R^1.prop}.1)  to
$f$ and  $f/\mu$ to get isomorphisms
$$
\q_X\stackrel{\tau_M^{-1}}{\longrightarrow}
 R^1(f/\mu)_*\q_{M/\mu}\stackrel{\pi^*}{\longrightarrow}
 R^1f_*\q_{M}\stackrel{\tau_{M/\mu}}{\longrightarrow} \q_X,
$$
whose composite is multiplication by $m(X)$.

Thus the map $\pi$ induces isomorphisms between the  spectral sequences
$$
H^i(X, R^jf_*\q_{M/\mu})\Rightarrow H^{i+j}(M/\mu,\q)
\qtq{and}
H^i(X, R^jf_*\q_M)\Rightarrow H^{i+j}(M,\q),
$$
where the maps
$$
\pi^*:H^i(X, R^jf_*\q_{M/\mu})\to H^i(X, R^jf_*\q_M)
$$
should be thought of as multiplication by $m(X)$.

Since $f/\mu: M/\mu\to X$ is a locally trivial circle bundle,
the edge homomorphisms
$$
 H^i(X, R^1(f/\mu)_*\q_{M/\mu})\to H^{i+2}(X,\q)
$$
are cup product with $c_1((M/\mu)/X)$.
Since $c_1(M/X)=\frac1{m(X)}c_1((M/\mu)/X)$
we see that the edge homomorphisms 
$$
\delta_i: H^i(X, R^1f_*\q_{M})\to H^{i+2}(X,\q)
$$
are cup product with $c_1(M/X)$.

Furthermore, 
$$
\tau(\pi^* H^0(X, R^1(f/\mu)_*\z_{M/\mu}))=m(X)\cdot H^0(X,\z)
$$
and so it agrees with $\tau(H^0(X, R^1f_*\z_M))$.
\end{say}

Thus we obtain the following:

\begin{cor}\label{chernclass=edgemap.cor}
Notation as above.
 The quotient map $\pi:M\to M/\mu$ (defined in (\ref{seif.maps.say}))
induces an isomorphism
$$
H^0(X, R^1f_*\z_{M})\cong \pi^* H^0(X, R^1(f/\mu)_*\z_{M/\mu}).
$$
Modulo torsion,  the edge homomorphisms
$$
\delta_i: H^i(X, R^1f_*\z_{M})\to H^{i+2}(X,\z)
$$
are identified with cup product with the Chern class  $c_1(M/X)$.
If $X$ is connected,  the image of
$$
\delta: H^0(X, R^1f_*\z_{M})\to H^{2}(X,\z)
$$
is generated by   $c_1(M/\mu)=m(X)c_1(M/X)$.\qed
\end{cor}

It is not clear to me how to describe the edge homomorphisms
on the torsion. Since $c_1(M/X)$ is not an integral class,
I do not even have a plausible guess.

Looking at the beginning of the Leray spectral sequence,
we obtain:

\begin{cor}\label{H^1=0.cond.cor} Notation and assumptions as above.
\begin{enumerate}
\item  $H^1(M,\q)=0$ iff $H^1(X,\q)=0$ and $c_1(M/X)\neq 0$.
\item If $H^1(M,\q)=0$ then $\dim H^2(M,\q)=\dim H^2(X,\q)-1$. \qed
\end{enumerate}
\end{cor}

 We see  that (\ref{R^1.prop}) describes
the sheaf $R^1f_*\z_M$ completely in terms of $(X,\Delta)$,
but  it is not always easy to compute its cohomologies 
based on this description. There are, however, some cases where
this is  quite straightforward.

\begin{prop}\label{R^1.compute.prop}
Let  $M$ be a differentiable manifold with a circle action
$\sigma:S^1\times M\to M$ with orientable stabilizer representations
and 
  $f:M\to (X,\sum (1-\frac1{m_i})D_i)$ 
the corresponding  Seifert bundle.
Set
$$
K:=\ker\bigl[ \z_X \to \sum_i \z_{D_i}/m_i\bigr].
$$
Then
\begin{enumerate}
\item there is an injection $\tau: R^1f_*\z_M\into K$ with quotient sheaf $Q$,
\item $\supp Q$ is the set of non--manifold points of $X$, and
\item $\dim \supp Q\leq \dim X-4$.
\end{enumerate}
\end{prop}

Proof. Pick a point $x\in X$ and let $m(x)$ denote the 
multiplicity of the Seifert fiber over $x$.
Pick a small neighborhood $x\in V_x$. Then
$H^0(V_x, R^1f_*\z_M)=m(x)\z$ by (\ref{R^1.prop}).

Let $C(x)$ be the product of those $m_i$ for which $x\in D_i$
and note that by (\ref{cod.4.prop}), $m_i$ and $m_j$ are relatively prime if 
$D_i\cap D_j\neq \emptyset$. 
Thus $H^0(V_x,K)=C(x)\z$.
By (\ref{cod.4.prop}), $C(x)$ divides $m(x)$ and
$X$ is a manifold at $x$
iff $m(x)=C(x)$.
\qed
\medskip

This allows us to compute some of the
cohomology groups of $R^1f_*\z_M$.

\begin{say}\label{bigcoh.comp} From (\ref{R^1.compute.prop}.3)
we conclude that $H^i(X, R^1f_*\z_M)=H^i(X,K)$ for $i\geq \dim X-2$
and we have  a long exact sequence computing $H^i(X,K)$.
Thus we get information on the 3 top cohomology groups
$H^i(X, R^1f_*\z_M)$. Set $\dim X=d$.

The top cohomology is the easiest:
$$
H^d(X, R^1f_*\z_M)\cong H^d(X,\z).
\eqno{(\ref{bigcoh.comp}.1)}
$$

For the next one, we have an exact sequence
$$
H^{d-2}(X,\z)\to \sum_i H^{d-2}(D_i,\z/m_i)\to H^{d-1}(X, R^1f_*\z_M)
\to H^{d-1}(X,\z)
$$
If $M$ is orientable, then the $D_i$ are orientable 
if $m_i\geq 3$, thus $H^{d-2}(D_i,\z/m_i)\cong \z/m_i$ for every $i$.
Thus we obtain:
\medskip

\ref{bigcoh.comp}.2 {\it Claim.} If $M$ is orientable 
and  $H^{d-1}(X,\z)=0$ then
$$
\begin{array}{rcl}
H^{d-1}(X, R^1f_*\z_M)
&=&\coker\bigl[H^{d-2}(X,\z)\to \sum_i H^{d-2}(D_i,\z/m_i)\bigr]\\
&=&\coker\bigl[H^{d-2}(X,\z)\to \sum_i \z/m_i\bigr].\qed
\end{array}
$$

The last relevant piece of the long exact sequence is
$$
H^{d-3}(X,\z)\to \sum_i H^{d-3}(D_i,\z/m_i)\to H^{d-2}(X, R^1f_*\z_M)
\to H^{d-2}(X,\z)
$$
We are especially interested in the torsion in
$H^{d-2}(X, R^1f_*\z_M)$. Most of it is coming from
$\sum_i H^{d-3}(D_i,\z/m_i)$, but 
$H^{d-3}(X,\z)$ and the torsion in $H^{d-2}(X,\z)$ influence it.
In general these are hard to control, but
we get the following:
\medskip

\ref{bigcoh.comp}.3 {\it Claim.} Assume that  $M$ is orientable, 
$H^{d-3}(X,\z)=0$ and $H^{d-2}(X,\z)$ is torsion free. Then
$$
H^{d-2}_{tors}(X, R^1f_*\z_M)\cong 
\sum_i H^{d-3}(D_i,\z/m_i).\qed
$$
\end{say}

\begin{say}[Proof of (\ref{main.thm}.1)] Let $L$ be a compact
5--manifold with a fixed point free circle action
with orientable stabilizer representations. 
Let 
$$
f:L\to \bigl(X, \sum \bigl(1-\tfrac1{m_i}\bigr)D_i\bigr)
$$
be the corresponding Seifert bundle.
By (\ref{dim5.nice.say}), $X$ has finitely many non--manifold points $X^s$
with complement $X^0=X\setminus X^s$.

Any abelian cover of $X^0$ gives an abelian cover
of $L\setminus L^s$, which then extends to an
abelian cover of $L$. Thus we conclude that
$H_1(L,\z)=0$ implies that 
$H_1(X^0,\z)=0$. By Lefschetz duality the latter gives that
$H^3(X,\z)=H^3(X,X^s,\z)\cong H_1(X^0,\z)=0$.
Furthermore, the torsion in
$H^2(X,\z)$ is isomorphic to the
torsion in $H_1(X,\z)$ hence again zero.

Let $b_1(D_i)$ denote 
 $\dim H_1(D_i,\z/2)$.
Thus $H_1(D_i,\z/m_i)=(\z/m_i)^{b_1(D_i)}$
since $D_i$ is orientable whenever $m_i\geq 3$.

Thus if $H_1(X^0,\z)=0$ then  the $E^2$-term of the 
 Leray spectral sequence
$$
H^i(X, R^jf_*\z_L)\Rightarrow H^{i+j}(L,\z)
$$
has the form
$$
\begin{array}{lcccl}
\z \quad & (\mbox{torsion}) & \quad\z^{k+1}+\sum_i (\z/m_i)^{b_1(D_i)}\quad & 
H^3(X,R^1f_*\z_L)\quad & \z\\[2ex]
\z & 0 & \z^{k+1} & 0 & \z.
\end{array}
$$

One can read off the cohomology of $L$ from this
spectral sequence. Let us start with a criterion
for the vanishing of $H_1(L,\z)$.

\begin{prop}\label{H_1=0.crit}
 Let $L$ be a compact, orientable
5--manifold with   a 
 Seifert bundle structure
$f:L\to \bigl(X, \sum \bigl(1-\tfrac1{m_i}\bigr)D_i\bigr)$
such that $H^3(X,\z)=0$.
\begin{enumerate}
\item There is a surjection
$$
H_1(L,\z)\onto H^3(X,R^1f_*\z_L)=
\coker\bigl[H^2(X,\z)\to  \sum_i H^2(D_i, \z/m_i)\bigr].
$$
\item Assume that  $H^3(X,R^1f_*\z_L)=0$ and 
 $X$ is smooth. Then the order of $H_1(L,\z)$ is $d$ iff
$$
c_1(L/\mu)=d\cdot(\mbox{primitive cohomology class})\in H^2(X,\z),
$$
(where primitive $:=$ not a nontrivial
multiple of any  cohomology class.)
\item Thus if $X$ is smooth then $H_1(L,\z)=0$ iff 
\begin{enumerate}
\item $H^2(X,\z)\to  \sum_i H^2(D_i, \z/m_i)$ is surjective, and
\item $c_1(L/\mu)\in H^2(X,\z)$ is primitive.
\end{enumerate}
\end{enumerate}
\end{prop}

Proof. By duality, $H_1(L,\z)\cong H^4(L,\z)$.
The spectral sequence shows that
$H^4(L,\z)\onto H^3(X,R^1f_*\z_L)$, hence 
 using (\ref{bigcoh.comp}.2) we obtain
the first claim.

By (\ref{chernclass=edgemap.cor}),  $c_1(L/\mu)\in H^2(X,\z)$
generates the image of the differential
$$
\delta_0:\z\cong H^0(X,R^1f_*\z_L)\to H^2(X,\z).
$$
 Thus if $c_1(L/\mu)=d\cdot \beta$
is not a primitive element, then $E^3_{0,2}$ contains $d$-torsion,
which survives in $H^2(L,\z)$, hence
the order of  $H_1(L,\z)$ is at least $d$.

Assume that $X$ is smooth and write
$c_1(L/\mu)=d\cdot \beta$ where $\beta$ is  primitive.
Since  cup product is a perfect pairing on  $H^2(X,\z)$,
there is an $\alpha\in H^2(X,\z)$ such that
$m(X)c_1(L)\cup \alpha=d$.
Thus $c_1(L)\cup m(X)\alpha=d$.
Since $m_i|m(X)$ for every $i$,
$$
m(X)H^2(X,\z)\subset \ker\bigl[H^2(X,\z)\to \sum_i H^2(D_i,\z/m_i)\bigr]
$$
and this kernel is the $\z^{k+1}$ summand of
$H^2(X,R^1f_*\z_L)$.
Hence $d=c_1(L)\cup m(X)\alpha$ is in the image of
$$
\delta_2: H^2(X,R^1f_*\z_L)\stackrel{c_1(L/X)\cup}{\longrightarrow} H^4(X,\z).
$$
Thus the order of $H_1(L,\z)$ also divides $d$.\qed

\medskip

\begin{cor}\label{H62.rest.onto}
 Let $L$ be a compact
5--manifold with  $H_1(L,\z)=0$, $\dim H_2(L,\q)=k$  and a 
 Seifert bundle structure
$f:L\to \bigl(X, \sum \bigl(1-\tfrac1{m_i}\bigr)D_i\bigr)$.
Then
$\#\{ i: p| m_i\}\leq k+1$. 
\end{cor}

Proof. Looking at  (\ref{H_1=0.crit}.1)
modulo $p$,
we obtain surjections
$$
\z^{k+1}\cong H^2(X,\z)\onto  \sum_i H^2(D_i, \z/m_i)
\onto  \sum_{i:p|m_i} \z/p.
$$
Thus $\#\{ i: p| m_i\}\leq k+1$. \qed
\medskip

From the spectral sequence we also get that
the torsion in $H^3(L,\z)$ is isomorphic to the
torsion in $H^3(X, R^1f_*\z_M))$ which in turn is
computed in  
(\ref{bigcoh.comp}.3).
Thus we obtain:

\begin{prop}\label{H_2.tors.onto}
 Let $L$ be a compact
5--manifold with  $H_1(L,\z)=0$ and a 
 Seifert bundle structure
$f:L\to \bigl(X, \sum \bigl(1-\tfrac1{m_i}\bigr)D_i\bigr)$.
Then there are isomorphisms
$$
\sum_i (\z/m_i)^{b_1(D_i)} 
\cong \sum_i H^1(D_i,\z/m_i)
\cong H^3_{tors}(L,\z).
$$
Dually, we can construct a basis of
$H_{2,tors}(L,\z)$ as follows:

Choose loops $\gamma_{ij}\subset D_i$, giving a basis
of $H_1(D_i,\z/m_i)$.
Then $\Gamma_{ij}:=f^{-1}(\gamma_{ij})\subset L$ is a 2--cycle
which is $m_i$-torsion and
$$
H_{2,tors}(L,\z)=\sum_{ij}(\z/m_i)[\Gamma_{ij}]. \qed
$$
\end{prop}

If $p^{a_i}$ is the largest $p$ power dividing $m_i$
then the $p$ part of $H_2(L,\z)$ is
$$
\sum_i (\z/p^{a_i})^{b_1(D_i)},
$$
and by (\ref{H62.rest.onto}) there are at most $k+1$ summands.
This proves (\ref{main.thm}.1).\qed
\end{say}

\section{Construction of Seifert bundles}

\begin{defn}
Let $X$ be a manifold and $D\subset X$ a codimension 
2 closed submanifold
with a tubular neighborhood $D\subset U\subset X$ and
 oriented normal bundle $N_D$. 
Thus we can view $N_D$ as a complex line bundle over $D$.
Choose an  identification
$j:U\cong N_D$.
 Composing with the bundle map $N_D\to D$
gives a retraction $\pi:U\to D$.
Thus $\pi^*N_D$ is a complex line bundle over $U$
and it has a section $s_D:u\to (u,j(u))$ which is
nowhere zero on $U\setminus D$. We can thus glue
$\pi^*N_D$ with the trivial complex line bundle on $X\setminus D$
to get a  complex line bundle $\o_X(D)$ with a section that vanishes
along $D$.

If $D=\sum c_iD_i$ is a formal integral linear
combination of 
codimension 2 closed submanifolds $D_i$ 
with oriented normal bundles $N_{D_i}$, then we define
$$
\o_X(D):=\bigotimes_i \o_X(D_i)^{\otimes c_i}.
$$
If $c_i\geq 0$ then $\o_X(D)$ has a natural section
which vanishes along $D_i$ with multiplicity $c_i$.

One can also define the cohomology class
of $D$ by $[D]:=c_1(\o_X(D))\in H^2(X,\z)$. 
The map $D\mapsto [D]$ is linear in the $c_i$, hence
it extends to rational linear combinations
giving $[D]\in H^2(X,\q)$. 
\end{defn}

\begin{thm}\label{seif.exist.thm} Let $X$ be a manifold and 
   $D_i\subset X$ codimension 2 closed submanifolds
with oriented normal bundles.
 Let  $1\leq b_i<m_i$ be  integers and $B$ a  complex line bundle on $X$.
Assume that
\begin{enumerate}
\item $(b_i,m_i)=1$  for every $i$,
\item  $(m_i,m_j)=1$ if $D_i\cap D_j\neq \emptyset$, and 
\item the $D_i$ intersect transversally.
\end{enumerate}
\noindent Then: 
\begin{enumerate}\setcounter{enumi}{3}
\item There is a   Seifert bundle  
$f:M=M(B,\sum\tfrac{b_i}{m_i}D_i)\to X$
such that
\begin{enumerate}
\item it has   orbit invariants
 $D_i,m_i,b_i$, and
\item $c_1(M/X)=c_1(B)+\sum\tfrac{b_i}{m_i}[D_i]$. 
\end{enumerate}
\item Every  Seifert bundle 
with orbit invariants
 $D_i,m_i,b_i$ is of the above form. 
\item The set of all  such Seifert bundles   forms a 
principal homogeneous
space under $H^2(X,\z)$ where the action  corresponds to
changing $B$.
\item The  properties (4.a--b) uniquely determine $M$ iff $H^2(X,\z)$ is
torsion free.
\end{enumerate}
\end{thm} 

Proof. 
Write $D=\sum\frac{b_i}{m_i}D_i$ and choose $m>0$ such that
every $m_i$ divides $m$.
In order to construct $M(B,\sum\tfrac{b_i}{m_i}D_i)\to X$
start with the rank 2 complex vector bundle
$$
h:E:= \o_X(mD)\otimes B^{\otimes m}+\o_X(\textstyle{\sum} D_i)\otimes B\to X.
$$
Since $m\sum D_i\geq mD$, the
natural section of $\o_X(m\sum D_i-mD)$ gives a  map
$$
\sigma:\o_X(mD)\otimes B^{\otimes m}\to \o_X(m\textstyle{\sum} D_i)
\otimes B^{\otimes m}.
$$
Define an auxiliary topological space  $N\subset E$ 
to be the set of all points
$$
\{(t,u,x): h(u)=h(t)=x\ \mbox{and}\ u^m=\sigma(t)\}.
$$
We see that $N$ is usually not normal, but we write down its normalization
$\bar N\to N$ explicitly, and we show that 
$\bar N\setminus(\mbox{zero section})$
is a Seifert $\c^*$-bundle whose unit circle bundle is
$M(B,\sum\tfrac{b_i}{m_i}D_i)$.

The key point is to get the local structure of $N$.

For  $x\in X$, one can choose an open  neighborhood
 in the form
$x\in U_x\cong \c^k+\r^{n-2k}$ with complex coordinates $x_1,\dots,x_k$
such that  $(x_i=0)$ are the components of $\sum D_i$ that pass through $x$.
After reindexing the $D_i$, we can assume that
$D_i=(x_i=0)$ with orbit invariants
$(m_i,b_i)$. 

Over this chart, we can write $\sigma(t)=t\prod x_i^{m(1-b_i/m_i)}$
thus $N$ is locally defined by the equation
$$
u^m=t\prod x_i^{m(1-b_i/m_i)}.
$$
Set $m_x=m_1\cdots m_k$, $m=m_xm'$ and define the $c_i$ by
the condition
$$
\sum c_ib_i\tfrac{m_x}{m_i}\equiv -1\mod m_x.
$$
This is solvable since
$\gcd( m_x/m_1,\dots,m_x/m_k)=1$.
In the notation of (\ref{cod.4.say}), 
$$
\begin{array}{l}
\ \ \ M_x:=\bigl(\c_v+
\c^k_{z_1,\dots,z_k}\bigr)/\tfrac1{m_x}
(1;-c_1\tfrac{m_x}{m_1},\dots,-c_k\tfrac{m_x}{m_k})
+\r^{n-2k}\\
f_x\downarrow\\ 
\ \ \ U_x:=\c^k_{z_1,\dots,z_k}/\tfrac1{m_x}
(c_1\tfrac{m_x}{m_1},\dots,c_k\tfrac{m_x}{m_k})
+\r^{n-2k}
\end{array}
$$
is a Seifert $\c^*$-bundle extended by the zero section. Observe that the 
 functions 
$$
s=v^{m_x}, x_i:=z_i^{m_i}\qtq{and} u:=v\prod z_i^{m_i-b_i}
$$
are invariant under the $\z/m_x$-action.
Set $t:=s^{m'}$. Then
$$
u^m=v^m\prod z_i^{m(m_i-b_i)}=t\prod x_i^{m(1-b_i/m_i)},
$$
hence we obtain a map $M_x\to N\cap h^{-1}(U_x)$ which 
gives the normalization. Thus the $M_x$ patch together
to $M\to N$ which is 
a Seifert $\c^*$-bundle extended by the zero section.
Furthermore, composing with the first projection
$pr_1:E\to \o_X(mD)\otimes B^{\otimes m}$ we get
$\pi:M\to \o_X(mD)\otimes B^{\otimes m}$ which is
the quotient of $M$ by $\z/m$.
Thus the first Chern class of $M/X$ is 
$$
c_1(M/X)=\tfrac1{m}(c_1(mD)+c_1(B^{\otimes m}))=c_1(D)+c_1(B).
$$
The rest follows from (\ref{seif.obstr}). \qed

\begin{cor}\label{H^2.rel.cor} Notation as in
(\ref{seif.exist.thm}). Set $E_i:=f^{-1}(D_i)$.
The normal bundle of  $E_i\subset M$ is orientable 
consistently with the normal bundle of $D_i$
and $f^*c_1(B)+\sum_i b_i[E_i]=0$.
\end{cor}

Proof. Consider the projection maps
$p:M\to \o_X(\sum D_i)\otimes B$ and
$h_2: \o_X(\sum D_i)\otimes B\to X$.
The pull back $h_2^*\o_X(\sum D_i)\otimes B$ has a tautological section
$U$ which vanishes only along the zero section.
In the local charts used in the proof of 
(\ref{seif.exist.thm}) this section is denoted by $u$.
From the formula $u=v\prod z_i^{m_i-b_i}$
we see that $p^*U$ vanishes along $E_i$ with
multiplicity $m_i-b_i$. Thus
$$
c_1(f^* \o_X(\sum D_i)\otimes B)=\sum (m_i-b_i)[E_i].
$$
Since $f^*[D_i]=m_i[E_i]$, this becomes
$f^*c_1(B)=-\sum b_i[E_i]$.\qed

\medskip

 Let $X$ be a topological space.
Continuous sections of $S^1\times X\to X$ form a sheaf,
 denoted by $S^1_X$. Its
cohomology groups are denoted by $H^i(X,S^1)$.

 Let $C^0_X$ denote the sheaf of continuous 
 functions. This sheaf is soft   and so it has no higher cohomologies
(cf.\ \cite[II.9]{bredon}).
Thus the 
 long exact cohomology sequence of
$$
0\to \z_X\to C^0_X\to S^1_X
\to 0
$$ 
shows that $H^i(X,S^1)\cong H^{i+1}(X,\z)$ for $i\geq 1$.

The following rather standard  result, closely related to
\cite[4.5]{HaSa91}, provides an 
approach to the global description of Seifert bundles.

\begin{prop} \label{seif.obstr}
 Let $X$ be a   topological  space and 
 $X=\cup U_i$  an open cover.  
Assume that over each $U_i$ we have a
 Seifert bundle $Y_i\to U_i$
and  there are $S^1$-equivariant     homeomorphisms
$\phi_{ij}:Y_j|_{U_{ij}}\cong Y_i|_{U_{ij}}$.

\begin{enumerate}
\item There is an obstruction element in the torsion subgroup
 $H^2_{tors}(X,S^1)\cong H^3_{tors}(X,\z)$
such that there is a global Seifert bundle $Y\to X$ compatible with these
local structures iff the obstruction element is zero.
\item The set of all such global Seifert bundles,  up to 
$S^1$-equivariant  homeomorphisms, 
 is either empty or forms a principal homogeneous
space under $H^1(X,S^1)\cong H^2(X,\z)$.
\item The action of $H^2(X,\z)$ on the Chern classes
is addition.
\end{enumerate}
\end{prop}

Proof. By (\ref{seif.maps.say}), 
the isomorphisms $\phi_{ij}$ can be changed
only to  $\alpha_{ij}\phi_{ij}$ for any
$\alpha_{ij}\in H^0(U_{ij},S^1)$.  These patchings define a 
global Seifert bundle iff
$$
\alpha_{ik}\phi_{ik}=\alpha_{ij}\phi_{ij}\alpha_{jk}\phi_{jk}
\qtq{for every $i,j,k$.}
$$
This is equivalent to
$$
\alpha_{ij}\alpha_{jk}\alpha_{ki}=(\phi_{ij}\phi_{jk}\phi_{ki})^{-1}
\qtq{for every $i,j,k$.}
\eqno{(\ref{seif.obstr}.3)}
$$
The products $(\phi_{ij}\phi_{jk}\phi_{ki})^{-1}\in H^0(U_{ijk}, S^1)$
satisfy the cocycle condition, and they define
an element of $H^2(X,S^1)$, called the obstruction.
One can find
$\{\alpha_{ij}\}$ satisfying
(\ref{seif.obstr}.3) iff the obstruction is zero.

Replacing the $Y_j$ by $Y_j/\mu$ as in (\ref{seif.maps.say}) changes the
isomorphisms over $U_i\cap U_j$ to $\phi_{ij}^M$, hence the
obstruction corresponding to the Seifert bundles 
$Y_j/\mu$ is the $M$th power of the original obstruction.

The quotients
 $Y_j/\mu$ are all $S^1$-bundles, and these can always be
globalized to the trivial $S^1$-bundle. Thus the obstruction is
 torsion.

Two choices  $\{\alpha_{ij}\}$ and $\{\alpha'_{ij}\}$
give isomorphic Seifert bundles iff there are
isomorphisms $\delta_i:Y_i\cong Y_i$
(viewed as elements of $H^0(U_i,S^1)$)
such that
$$
\alpha'_{ij}\alpha^{-1}_{ij}=\delta_i\delta^{-1}_j|_{U_{ij}}.
$$
Thus  $\{\alpha'_{ij}\alpha^{-1}_{ij}\}$
corresponds to a class in $H^1(X,S^1)$
and we also get (3).\qed
\medskip

By (\ref{seif.exist.thm}) 
 the obstruction vanishes if $X$  and the $D_i$ are
orientable,  but
there are even complex orbifold examples where the obstruction is nonzero
\cite[Exmp.35]{ko-seif}.

\section{The second Stiefel--Whitney class}

We start the computation of the  second Stiefel--Whitney class
of a Seifert bundle by two key examples.

\begin{exmp}\label{loops.exmps}
 Given $1\leq b <m$ such that $(m,b)=1$
consider the map
$$
f: S^1_s\times \c_z\times \c^*_y\mapsto \c_x\times\c^*_y
\qtq{given by} (s,z,y)\mapsto (s^bz^m, y),
$$ 
It has a Seifert  bundle structure
where the $S^1_t$-action is given by
$$
(t)\times (s,z,y)\mapsto (t^{-m}s,t^bz, y).
$$
For any $y_0$, we have  $f^{-1}(0,y_0)=(*,0,y_0)$
and $\stab_{(*,0,y_0)}=\z/m$. A transverse 
 slice at $(s_0,0,y_0)$ is given by
$(s_0,*,y_0)$ and the stabilizer representation
is $z\mapsto \epsilon^{b}z$ where $\epsilon=e^{2\pi i/m}$.

Note that the tangent bundle of $S^1_s\times \c_z\times \c^*_y$
is parallelizable, hence its Stiefel--Whitney classes are all zero.

Since $H^2(\c_x\times\c^*_y,\z)=0$, we conclude from (\ref{seif.obstr})
that the
above examples exhaust all possible Seifert bundles
over the orbifold
$(\c_x\times\c^*_y, (1-\frac1{m})D)$,
where $D=\{0\}\times \c^*_y$.
\end{exmp}

\begin{exmp}\label{loops.nonor.exmps}
 For  $(m,b)=(2,1)$ 
there is also a version where the branch divisor
is not orientable. Indeed, consider the 
orientation preserving involutions
$$
(s,z,y)\mapsto (s, \bar s\bar z, -1/\bar y)
\qtq{and}
(x,y)\mapsto (\bar x, -1/\bar y).
$$
These commute with the map $f$ and the $S^1_t$-action in (\ref{loops.exmps}).
Thus we get a Seifert bundle structure on  the quotient.

Since $H^2((\c_x\times\c^*_y)/(\z/2),\z)=0$, as
above we see that 
this is the only Seifert bundle
over the orbifold
$((\c_x\times\c^*_y)/(\z/2), (1-\frac1{2})D)$, where
$D=(\{0\}\times \c^*_y/(\z/2))$ is not orientable.

It is easy to see that the second Stiefel--Whitney class
is nonzero on the 2--cycle
$(z=0, |y|=1)/(\z/2)$. 
\end{exmp}

The next result proves  (\ref{main.thm}.2)
and also gives more information about  the
invariant $i(L)$ defined in 
(\ref{i(L).defn}).

\begin{prop}\label{w_2.ontors.prop} Let $L$ be a
5--manifold with $H_1(L,\z)=0$ having a Seifert bundle structure
$f:L\to \bigl(X, \sum \bigl(1-\tfrac1{m_i}\bigr)D_i\bigr)$.
Then
\begin{enumerate}
\item $i(L)\in \{0,1,\infty\}$, and
\item $i(L)=1$ iff at least one of the $D_i$ is
nonorientable.
\end{enumerate}
\end{prop}

Proof.
Fix a $D_i$  and
choose a loop $\gamma\subset D_i^0$.
By (\ref{H_2.tors.onto}), $\Gamma:=f^{-1}(\gamma)\subset L$ is a 2--cycle
which is $m_i$-torsion in $H_2(L,\z)$ 
and these cycles generate the torsion subgroup of $H_2(L,\z)$.

Let $\gamma\subset V\subset X^0$ be a tubular neighborhood.
If $D_i$ is orientable along $\gamma$, then
the pair $(V, D_i\cap V)$ is diffeomorphic to 
$(\c_x\times\c^*_y, \{0\}\times \c^*_y)$, thus the restriction
of $f:L\to X$ to $V$ is diffeomorphic to one of the
Seifert bundles enumerated in (\ref{loops.exmps}).
Therefore $w_2(L)\cap [\Gamma]=0$.
Since $m_i\geq 3$ implies that $D_i$ is orientable
(\ref{norm.bund.orient}),
we get the first claim.

If every $D_i$ is orientable then we get that $w_2$
is zero on all the torsion, hence
$i(L)\in \{0,\infty\}$. Conversely, if 
$D_i$ is nonorientable along $\gamma$, then
the pair $(V, D_i\cap V)$ is diffeomorphic to 
the one in (\ref{loops.nonor.exmps}), hence
$w_2(L)\cap [\Gamma]\neq 0$ and so $i(L)=1$.
\qed
\medskip

We also need the following formula for
the second Stiefel--Whitney class of a Seifert bundle.
It is the topological version of the formula for
the first Chern class for holomorphic Seifert bundles given in
\cite{fl-za} and
\cite[Cor.41]{ko-seif}. 

\begin{lem} 
Let $M$ be a
manifold with a fixed point free circle action
with orientable stabilizer representations 
and 
$f:M\to \bigl(X, \sum \bigl(1-\tfrac1{m_i}\bigr)D_i\bigr)$
the corresponding Seifert bundle, $X$ smooth.
Set $E_i:=f^{-1}(D_i)$. Then 
$$
w_2(M)=f^*w_2(X)+\sum_i (m_i-1)[E_i].
$$
\end{lem}

Proof. We factor $f$ as the composite of $\pi:M\to M/\mu$
and of the projection $f/\mu:M/\mu\to X$. Since
$M/\mu\to X$ is a circle bundle, $T_{M/\mu}=(f/\mu)^*T_X+
(\mbox{trivial bundle})$, thus $w_i(M/\mu)=f^*w_i(X)$ for every $i$.
Note that $\pi:M\to M/\mu$ is a branched covering which ramifies along the
subspaces $E_i$ and the ramification order is $m_i$.
Thus we need to show that 
$w_2(M)=f^*w_2(M/\mu)+\sum_i (m_i-1)[E_i]$.
For complex manifolds and for $c_1$ instead of $w_2$
this is the Hurwitz formula. 

For $\dim M\geq 5$, 
we can represent any homology class $H_2(M,\z/2)$ by an embedded
surface $S\into M$. We may also assume that $S$ is transversal
to $\cup E_i$ and near each point $E_i\cap S$
the induced map $\pi:S\to \pi(S)$ is  a branched cover of
degree $m_i$. In a neighborhood of $\pi(S)$ we can write
$T_{M/\mu}=E+N$ where $N$ is a trivial bundle, $\rank E=2$ 
and near each point $\pi(E_i\cap S)$ the subbundles
$T_{\pi(S)}\subset   T_{M/\mu}$ and $E \subset   T_{M/\mu}$
agree. Thus
$\pi^*w_2(M/\mu)\cap [S]=
w_2(\pi^*E)\cap S$.
Correspondingly one can write $T_M|_S=E'+\pi^*N$ 
where $E'$ is a rank 2 bundle with 
induced tangent map $E'\subset \pi^*E$ whose quotient
is supported at the points $E_i\cap S$ and has length $m_i-1$ there.\qed

\medskip

This formula becomes easier to use if we
 combine it with (\ref{H^2.rel.cor})
which says that $f^*c_1(B)+\sum b_i[E_i]=0$.
Thus
$$
\begin{array}{rcl}
w_2(M)&=&f^*w_2(X)+\sum_i (m_i-1)[E_i]+f^*c_1(B)+\sum b_i[E_i]\\
&=& f^*w_2(X)+\sum_i (m_i-1+b_i)[E_i]+f^*c_1(B).
\end{array}
$$
If $m_i$ is even then $b_i$ is odd so $m_i-1+b_i$ is even.
If $m_i$ is odd then $m_i-1$ is even, so we can rewrite this as
$$
w_2(M)=f^*w_2(X)+\sum_{i: m_i\ {\rm odd}} b_i[E_i]+f^*c_1(B).
$$
Note that in integral cohomology, $m_i[E_i]=f^*[D_i]$, hence,
in $H^2(M,\z/2)$, we get that  $[E_i]=f^*[D_i]$ 
if $m_i$ is odd and  $f^*[D_i]=0$ if $m_i$ is even.
Thus we can rewrite our formula as follows.

\begin{cor} \label{w_2=pullback}
Let 
$M=M(B,\sum \frac{b_i}{m_i})\to X$
be a Seifert bundle as in (\ref{seif.exist.thm}),
 $X$ smooth and $X,D_i$ orientable. Then
$$
w_2(M)=f^*\bigl(w_2(X)+\sum_i b_i[D_i]\bigr)+f^*c_1(B).\qed
$$
\end{cor}

\begin{say}[Proof of (\ref{main.thm}.3)]
Let  $f:L\to (X,\sum (1-\frac1{m_i})D_i)$,
be a  Seifert bundle such that $H_1(L,\z)=0$.
 If $i(L)=\infty$ then
every $D_i$ is orientable by (\ref{w_2.ontors.prop}).

By (\ref{dim5.nice.say}), there is a finite set $X^s\subset X$ such that
$X^0:=X\setminus X^s$ is a manifold.
Set $L^0:=f^{-1}(X^0)$. Since $L\setminus L^0$ has codimension
$4$, we see that $H^i(L,\z/2)=H^i(L^0,\z/2)$ for $i\leq 2$.

By (\ref{H_2.tors.onto}),  if $c(2^i)\neq 0$ for $k+1$ values of $i$ 
then $\geq k+1$ of the $m_j$ are even. Let these be
$D_0,\dots,D_k$. Since
$H^2(X,\z)\to \sum_i H^2(D_i,\z/m_i)$ is surjective,
we conclude that
$$
H^2(X,\z/2)\to \sum_{i=0}^k H^2(D_i,\z/2)\qtq{is surjective.}
$$
The two sides have the same rank, so we have an isomorphism.
This implies that $D_0,\dots,D_k$ form a basis of
$H_2(X,\z/2)=H_2(X,X^s,\z/2)$. 

By Lefschetz duality,
$\o_{X^0}(D_0),\dots, \o_{X^0}(D_k)$
(or rather their Chern classes) form a basis of $H^2(X^0,\z/2)$.
Since $f^*\o_{X^0}(D_j)=\o_{L^0}(m_jE_j)$, 
and $m_0,\dots,m_k$ are even, we conclude
that the pull back map
$f^*: H^2(X^0,\z/2)\to  H^2(L^0,\z/2)$ is zero.

On the other hand,
 (\ref{w_2=pullback}) applies and
$w_2(L^0)$ is the pull back of a
cohomology class from $X^0$. Thus $w_2(L)=w_2(L^0)=0$.
\qed
\end{say}

\section{Seifert bundles over $\c\p^2$}

In this section we construct  examples
of Seifert bundles $f:L\to (X,\Delta)$. 
The base $(X,\Delta)$ is constructed as a connected
sum of pairs $(\c\p^2,\sum (1-\frac1{m_i})D_i)$,
where the attaching does not involve the $D_i$.
Once the base $(X,\Delta)$ is chosen, we can vary
the complex line bundle $B$ to obtain various
values of the invariant $i(L)$.

We already have a good understanding of the cohomology
of $L$, the key additional step is to control
the fundamental group as well.
This is achieved by the following simple lemma.

\begin{lem}\label{fg.nilp.lem}  Let $f:L\to (X,\Delta)$
be a Seifert bundle, $X$ a manifold.
 If $\pi_1(X\setminus\Delta)$ is solvable then
so is $\pi_1(L)$.
If this holds then $\pi_1(L)=1$ iff $H_1(L,\z)=0$.
\end{lem}

Proof. Since $f:L\setminus f^{-1}(\Delta)\to X\setminus\Delta$
is a circle bundle, there is an exact sequence
$$
\pi_1(S^1)\to \pi_1(L\setminus f^{-1}(\Delta))\to
 \pi_1(X\setminus\Delta)\to 1.
$$
Since $f^{-1}(\Delta)\subset L$ has codimension 2,
there is a surjection 
$\pi_1(L\setminus f^{-1}(\Delta))\to \pi_1(L)$.\qed

\begin{rem} \label{fund.gr.sequence}
Although we do not need it, it is 
worthwhile to note that there is an exact sequence for
$\pi_1(L)$ itself.

 Let $(X,\sum_i (1-\frac1{m_i})D_i)$ be an orbifold
and  $X^0\subset X$  the smooth locus of $X$.
The {\it orbifold fundamental group} $\pi^{orb}_1(X,\Delta)$
is the fundamental group of $X^0\setminus\supp\Delta$
modulo the relations: if $\gamma$ is any 
 small loop around $D_i$ 
 then $\gamma^{m_i}=1$ \cite{thur}.

Note that $\pi^{orb}_1(X,\emptyset)$ may be different from
$\pi_1(X)$ if $X$ is not a manifold.

The abelianization of $\pi^{orb}_1(X^0,\Delta)$,
denoted by  $H^{orb}_1(X^0,\Delta)$, is called the
{\it abelian orbifold fundamental group.}
(The higher orbifold homotopy and homology groups
are defined in \cite{Hae84}.)

A straightforward generalization of the
computation of the fundamental group of 3--dimensional Seifert 
bundles (see \cite{seif} or \cite[5.7]{HaSa91})
gives the 
 exact sequence
$$
\pi_1(S^1)\to \pi_1(L)\to \pi^{orb}_1(X,\Delta)\to 1.
$$
\end{rem}

The next lemma gives a large collection of pairs
$(\c\p^2, \sum D_i)$ to work with.

\begin{lem}\label{D_i.in.CP2} Let $D_1,\dots,D_s$ be compact surfaces.
Then there are embeddings
$D_i\subset \c\p^2$ such that
\begin{enumerate}
\item the $D_i$ intersect transversally,
\item if $D_i$ is orientable, then its homology class
$[D_i]$ is a  generator of $H_2(\c\p^2,\z)$,
\item if $D_i$ is nonorientable, then 
$[D_i]$ is a  generator of $H_2(\c\p^2,\z/2)$,
\item  $\pi_1(\c\p^2\setminus(D_1\cup\cdots\cup D_s))$ 
is abelian.
\end{enumerate}
\end{lem}

Proof. Let us start with $\c^2$ and for each surface
with $b_1(D_i)$ even pick a complex line $L^0_i$  in general position
and for each surface
with $b_1(D_i)$ odd pick a non complex real affine 2--plane 
$L^0_i\subset \c^2$ in general position.
Correspondingly, in $\c\p^2$ we get
embedded copies $L_1,\dots,L_s$ of $\c\p^1$ and of $\r\p^2$
which intersect transversally and
satisfy the conditions (\ref{D_i.in.CP2}.2--3).

A classical lemma \cite[p.317]{zar}
states (in the case of complex lines)
that $\pi_1(\c^2\setminus(L^0_1\cup\cdots\cup L^0_s))$ is abelian,
 but the proof
applies to real 2--planes in $\r^4$ as above
as long as all intersections are  transverse and  any two of 
the planes  do
intersect.
This implies that 
  $\pi_1(\c\p^2\setminus(L_1\cup\cdots\cup L_s))$ is abelian.

Next we aim to attach handles to the $L_i$
without changing the fundamental group of the complement.
The key part is the following local computation.

\begin{say}[Attaching handles]\label{handles.exmp}
 Take $\r^4$ with coordinates $(x,y,z,t)$.
A pl--embedded copy of $\r^2$ is given by the union of
the two half planes
$H_1:=(x\leq 0, y,0,0)$ and $H_2:=(0,y,0,t> 0)$. 
Make 2 holes in $H_1$ and attach a handle $S^1\times [-1,1]$ inside
$\r^3\cong (t=0)$ to $H_1$. Depending on how this is done, the resulting
surface can be orientable or nonorientable.
Together with $H_2$, we obtain a pl--embedded surface
$D\subset \r^4$.

We claim that $\pi_1(\r^4\setminus D)\cong \z$,
generated by a loop around $H_2\subset (t> 0)$.

Indeed, any loop in $\r^4\setminus D$ can be made
transversal to the hyperplane $(t=0)$.
Since $(t=0)\setminus D$ is connected,
we can assume that the loop intersects
$(t=0)$ only at points where $x>0$.

Since $D$ is disjoint from the half space $(t<0)$,
any part of the loop in this half space can be contracted and
then pushed above $(t=0)$. Thus we homotoped the
loop to the upper half space $(t>0)$
and $\pi_1((t>0)\setminus H_2)\cong \z$.
\end{say}

To get the final embeddings,
pick  points $p_{i}\in L_i$ not on any of the other $L_j$  and
 disjoint  neighborhoods $p_{i}\in U_{i}\sim \r^4$
such that $L_i\cap U_{i}\into U_{i}$ is a linearly embedded $\r^2$.
We can attach handles to the $L_i$ to get
$D_i\subset \c\p^2$. In doing this, we have not created new
intersections and  the homology class
of $D_i$ is the same as the homology class of the $L_i$.
(With the caveat that if $L_i$ is orientable but $D_i$ is not
then we claim only a modulo 2 equality.)

Finally we need to show that
 $\pi_1(\c\p^2\setminus(D_1\cup\cdots\cup D_s))$ is abelian.
We prove that it is isomorphic to
 $\pi_1(\c\p^2\setminus(L_1\cup\cdots\cup L_s))$, which
we already know to be  abelian

We can compute 
both of these fundamental groups using van Kampen's theorem
from the $\pi_1(U_i\setminus L_i)$
(resp.\ $\pi_1(U_i\setminus D_i)$)  and
$$
\c\p^2\setminus(U_1\cup\cdots\cup U_s)\setminus(L_1\cup\cdots\cup L_s)
=\c\p^2\setminus(U_1\cup\cdots\cup U_s)\setminus(D_1\cup\cdots\cup D_s).
$$
Thus it is enough to prove that both are computed the same way,
which amounts to proving that 
$\pi_1(U_i\setminus L_i)=\pi_1(U_i\setminus D_i)$.
This was already done in (\ref{handles.exmp}).\qed

\medskip

{\it Note.} If the $D_i$ are all orientable,
although $D_i$ and its normal bundle both have a complex
structure, these can not be made compatible with the
(almost) complex
structure of $\c\p^2$. Indeed, 
with the usual  complex
structure of $\c\p^2$, 
$c_1(T_{\c\p^2}|_{D_i})=3$
and using the above complex
structures on $D_i$ and its normal bundle
gives 
$$
c_1(T_{\c\p^2}|_{D_i})=c_1(T_{D_i})+c_1(N_{D_i})=2-2g_i+1=3-2g_i.
$$
\medskip

\begin{const}[Seifert bundles with $i(L)\in \{0,\infty\}$]
\label{seif.0.inf.constr}
 Assume  that we have a natural number $k$ and a sequence of
natural numbers $c(p^i)$ such that 
\begin{enumerate}
\item the $c(p^i)$ are all even, and
\item for every prime $p$,  at most
$k+1$ of the $c(p^i)$ are nonzero.
\end{enumerate}

We want to construct
 Seifert bundles $f:L\to (X,\sum (1-\frac1{m_i})D_i)$ with 
$i(L)\in \{0,\infty\}$. By
(\ref{w_2.ontors.prop}.2) this means that every $D_i$ is orientable.

For every prime $p$ arrange the $p^i$ satisfying $c(p^i)\neq 0$
in increasing order, and put 1's at the beginning to get 
a sequence of length
$k+1$
$$
1=p^0,\dots, p^0< p^{i(1)}< \cdots< p^{i(h)}.
$$
Let $m_{pj}$ be the $j$th element of this sequence
and set $g_{pj}=\frac12c(m_{pj})$ where we set
$c(1)=0$. 
If $c(2^i)=0$ for every $i$ then 
set $m_{2j}=1$ for $j<k$, $m_{2k}=2$ and $g_{2j}=0$ for every $j$.

For each $j=0,\dots,k$,
pick a copy $\c\p^2_j$ of $\c\p^2$ and as in
(\ref{D_i.in.CP2})   construct
 oriented surfaces $D_{pj}\subset \c\p^2_j$ of genus $g_{pj}$.
Let $X$ be the connected sum of the $\c\p^2_j$
with $D_{pj}\subset X$. Thus
 $D_{pj}\cap D_{p'j'}=\emptyset$ iff $j\neq j'$.

Let $H_j\in H^2(X,\z)$ denote the cohomology class of a line on
$\c\p^2_j$. Then $H_0,\dots,H_k$ is a free basis
of $H^2(X,\z)$.

The restriction map
$H^2(X,\z)\to \sum_{pj} H^2(D_{pj},\z/m_{pj})$
can be written as a sum of the individual
 restriction maps
$$
\z\cong H^2(\c\p^2_j,\z)\to \sum_{p} H^2(D_{pj},\z/m_{pj})\cong
\sum_p \z/m_{pj}\cong \z/\bigl(\prod_p m_{pj}\bigr),
\eqno{(\ref{seif.0.inf.constr}.1)}
$$
where the last isomorphism holds since $m_{pj}$ and $m_{p'j}$ are
relatively prime for $p\neq p'$.
Thus the first condition of
(\ref{H_1=0.crit}.3) is satisfied.

On $\c\p^2_k$ we have surfaces $D_{pk}$ and assigned
multiplicities $m_{pk}$.
Let $m(X)=\prod_p m_{pk}$.
Since the $m_{pk}$ are relatively prime to each other
and their product is $m(X)$, we can find
$1\leq b_{pk}<m_{pk}$ such that $(b_{pk},m_{pk})=1$ and 
$$
\sum_p \frac{b_{pk}}{m_{pk}}\equiv \frac1{m(X)}\mod 1.
$$
Set $b_{pj}=1$ for $j<k$; these values are unimportant for us.
 
By (\ref{seif.exist.thm}), 
we can fix $h_k$ such that 
for any $h_0,\dots,h_{k-1}$ there is a Seifert bundle
$$
L(h_0,\dots,h_{k-1}):=
L\bigl(\o_X(\sum_j h_jH_j),\sum_{p,j} \frac{b_{pj}}{m_{pj}}D_{pj}\bigr)\to X
\eqno{(\ref{seif.0.inf.constr}.2)}
$$
and  $\gamma_j\in \bigl( \prod_p m_{pj}\bigr)^{-1}\z$ such that 
$$
c_1(L(h_0,\dots,h_{k-1}))=\frac1{m(X)}H_k+\sum_{j=0}^{k-1} \gamma_j H_j.
\eqno{(\ref{seif.0.inf.constr}.3)}
$$
Since $m_{pk}$ contains $p$ with the largest
exponent,  we see that $m_{pj}|m(X)$ for every $j$
and $m(X)/m_{pj}$ is even for every $j<k$.
Thus 
$$
\begin{array}{rcl}
c_1(L(h_0,\dots,h_{k-1})/\mu)&=&
m(X)c_1(L(h_0,\dots,h_{k-1}))\\
&\in & H_k+ 2\langle H_0,\dots,H_{k-1}\rangle.
\end{array}
\eqno{(\ref{seif.0.inf.constr}.4)}
$$
This implies that $c_1(L(h_0,\dots,h_{k-1})/\mu)$ 
is a primitive vector
and so the second condition of  (\ref{H_1=0.crit}.3)
also holds,  implying that  $H_1(L(h_0,\dots,h_{k-1}),\z)=0$. From
(\ref{D_i.in.CP2}) and (\ref{fg.nilp.lem}) we conclude that
$L(h_0,\dots,h_{k-1})$ is simply connected.
By (\ref{H_2.tors.onto}), 
$$
H_2(L(h_0,\dots,h_{k-1}),\z)=\z^k+\sum_{p,i} \bigl(\z/p^i\bigr)^{c(p^i)}.
$$
So far the $h_0,\dots h_{k-1}$ played no visible role in the construction.
Now we aim to choose these appropriately to control 
$w_2(L(h_0,\dots,h_{k-1}))$.

We have a pull back map $f^*:H^2(X,\z/2)\to H^2(L(h_0,\dots,h_{k-1}),\z/2)$
and by (\ref{w_2=pullback})
 there is an $\eta\in H^2(X,\z/2)$
such that 
$$
w_2\bigl(L(h_0,\dots,h_{k-1})\bigr)
=f^*\bigl(\eta+\sum_{i=0}^{k-1}h_jH_j\bigr).
\eqno{(\ref{seif.0.inf.constr}.5)}
$$
By choosing $h_0,\dots h_{k-1}$ suitably, we can thus assume that
$$
w_2\bigl(L(h_0,\dots,h_{k-1})\bigr)
=f^*(cH_k)\qtq{for some $c$.}
\eqno{(\ref{seif.0.inf.constr}.6)}
$$
 We know from (\ref{chernclass=edgemap.cor}) that
$c_1(L(h_0,\dots,h_{k-1})/\mu)$
 is in the kernel of $f^*:H^2(X,\z)\to H^2(L(h_0,\dots,h_{k-1}),\z)$,
hence its mod 2 reduction, which is $H_k$ by 
(\ref{seif.0.inf.constr}.4),
is in the kernel of 
$f^*:H^2(X,\z/2)\to H^2(L(h_0,\dots,h_{k-1}),\z/2)$.
Thus with the choices of (\ref{seif.0.inf.constr}.6),
we have a trivial second Stiefel--Whitney class.
\end{const}

We have completed the first part of the existence 
theorem (\ref{main.thm.2}):

\begin{cor} Let $c(p^i)$ be even natural numbers
(only finitely many nonzero)
satisfying (\ref{main.thm}.1) for some $k$. Then there are Seifert bundles 
$f:L\to ((k+1)\#\c\p^2,\Delta)$ 
such that $L$ is simply connected, $w_2(L)=0$ (equivalently,  $i(L)=0$) and
$$
H_2(L,\z)=\z^k+\sum_{p,i} \bigl(\z/p^i\bigr)^{c(p^i)}.\qed
$$
\end{cor}

Getting examples with $i(L)=\infty$ is very similar.
First, by changing the choice of $h_0$ in 
(\ref{seif.0.inf.constr}.6), we can achieve that
$$
w_2\bigl(L(h_0,\dots,h_{k-1})\bigr)
=f^*(H_0+cH_k).
\eqno{(\ref{seif.0.inf.constr}.7)}
$$
(The notation does not show, but we have to keep in mind
that $f^*$ also depends on $h_0,\dots,h_{k-1}$.) 
We know that $H_0+cH_k$ is not in the kernel of the
pull back map between integral cohomology groups.
 If we can show that  $H_0+cH_k$ is not in the kernel of the
pull back map between $\z/2$--cohomology groups, then we get
the desired examples with $i(L)=\infty$.

By assumption, at most $k$ of the $c(2^i)$ 
are nonzero, thus
there is an $i\geq 1$ such that every $m_{pj}$ is odd for
$j<i$ but $m_{2i}$ is even.

Let $Y$ be the connected sum of the first $i+1$ copies of
$\c\p_j$ and $g:M\to (Y,\Delta)$
the corresponding Seifert bundle.
By the above parity considerations,
$$
c_1(M/\mu)\in 2\langle H_0,\dots,H_{i-1}\rangle+\z H_i.
$$
Only one of the $m_{pj}$ with $j\leq i$ is even,
hence from the exact sequence given by (\ref{R^1.compute.prop})
$$
H^0(Y,\z)\to \sum_{j\leq i}H^0(D_{pj}, \z/m_{pj})\to H^1(Y,R^1g_*\z_M)\to
H^1(Y,\z)=0
$$
we conclude that $H^1(Y,R^1g_*\z_M)$ has odd order.
Thus the 
pull back map between integral cohomology groups sits in an exact
sequence
$$
\z\stackrel{c_1(M/\mu)}{\longrightarrow}
 H^2(Y,\z)\stackrel{g^*}{\to} H^2(M,\z)\to (\mbox{odd order group}).
$$
This implies that modulo 2 we still get an injection
$$
g^*:\langle H_0,\dots,H_{i-1}\rangle\into H^2(M,\z/2).
\eqno{(\ref{seif.0.inf.constr}.8)}
$$
Let $y\in Y$ be a point not on any of the $D_{pj}$
and set $Y^0:=Y\setminus \{y\}$ and $M^0:=M\setminus g^{-1}(y)$.
Then
$H^2(M,\z/2)=H^2(M^0,\z/2)$ and 
we can think of $M^0$ as an open subset of $L$. Thus
(\ref{seif.0.inf.constr}.8)
implies that 
$$
f^*:\langle H_0,\dots,H_{i-1}\rangle\into H^2(L,\z/2)
\qtq{is an injection. }
\eqno{(\ref{seif.0.inf.constr}.9)}
$$

Since $f^*(H_k)=0\in H^2(L,\z/2)$ by
(\ref{seif.0.inf.constr}.4), we conclude
that $f^*(H_0+cH_k)\in H^2(L,\z/2)$ is nonzero,
giving the second existence result:

\begin{cor} Let $c(p^i)$ be even natural numbers
(only finitely many nonzero)
satisfying (\ref{main.thm}.1) and (\ref{main.thm}.3)
 for some $k\geq 1$. Then there are Seifert bundles 
$f:L\to ((k+1)\#\c\p^2,\Delta)$ 
such that $L$ is simply connected, $i(L)=\infty$ and
$$
H_2(L,\z)=\z^k+\sum_{p,i} \bigl(\z/p^i\bigr)^{c(p^i)}.\qed
$$
\end{cor}

\begin{const}[Seifert bundles with $i(L)=1$]
\label{seif.1.constr}
 Assume  that we have a natural number $k$ and a sequence of
natural numbers $c(p^i)$ such that 
\begin{enumerate}
\item the $c(p^i)$ are all even except possibly $c(2)\geq 1$, and
\item for every prime $p$,  at most
$k+1$ of the $c(p^i)$ are nonzero.
\end{enumerate}

We want to construct
a Seifert bundle $f:L\to (X,\sum (1-\frac1{m_i})D_i)$ with 
$i(L)\in 1$. By
(\ref{w_2.ontors.prop}.2) 
we can assure this by choosing
 one of the  $D_i$  to be nonorientable.

The construction follows (\ref{seif.0.inf.constr}), with two changes:
\begin{enumerate}
\item We do not have to compute $w_2(L)$, since
$i(L)=1$ is guaranteed by having a nonorientable
surface $D_{2j}$.
\item  There is a unique $j$ with $m_{2j}=2$ and then we want
$D_{2j}\subset \c\p^2_j$ to be nonorientable with 
$H^1(D_{2j},\z/2)\cong (\z/2)^{c(2)}$. We can find such
$D_{2j}$ by (\ref{D_i.in.CP2}), but the Seifert bundle over
$\c\p^2_j$ has to be constructed by hand since the
 existence result (\ref{seif.exist.thm}) works only for
orientable surfaces.
\end{enumerate}

Let $N\subset \c\p^2$ be  a compact nonorientable surface
which is nonzero in $H_2(\c\p^2,\z/2)$. From the sequence
$$
H_2(N,\z/2)\stackrel{\cong}{\to} H_2(\c\p^2,\z/2)
\to H_2(\c\p^2,N,\z/2)
$$
and Lefschetz duality we conclude that
$H^2(\c\p^2,\z/2)
\to H^2(\c\p^2\setminus N,\z/2)$ is the zero map, thus
$$
\im\bigl[H^2(\c\p^2,\z)
\to H^2(\c\p^2\setminus N,\z)\bigr]\subset 2\cdot H^2(\c\p^2\setminus N,\z).
$$
By (\ref{seif.exist.thm}), there is a Seifert bundle
$$
g:M\to \bigl(\c\p^2, \sum_{p\geq 3}(1-\tfrac1{m_{pj}})D_{pj}\bigr)
$$
such that $c_1(M/\mu)\in H^2(\c\p^2,\z)$ is the generator.
Set $N=D_{2j}\subset \c\p^2$
and consider the restriction
$$
g:M\setminus g^{-1}(N)\to 
\bigl(\c\p^2\setminus N, \sum_{p\geq 3}(1-\tfrac1{m_{pj}})D_{pj}\bigr).
$$
As noted above
$$
c_1(M\setminus g^{-1}(N)/\mu)\in 
H^2(\c\p^2\setminus N,\z)
$$
is not primitive, but twice a generator. By
(\ref{H_1=0.crit}) this means that $\pi_1(M\setminus g^{-1}(N))=\z/2$
 and we have a ramified
 double cover
$$
g':M'\to M\to 
\bigl(\c\p^2, 
(1-\tfrac12)D_{2j}+\sum_{p\geq 3}(1-\tfrac1{m_{pj}})D_{pj}\bigr).
$$
This gives the required Seifert bundle with a
nonorientable $D_i$.
The rest of (\ref{seif.0.inf.constr}) works without changes and we obtain the
final existence result:
\end{const}

\begin{cor} Let $c(p^i)$ be even natural numbers for $p^i\geq 3$ 
(only finitely many nonzero)
satisfying (\ref{main.thm}.1) 
 for some $k$ and $c(2)\geq 1$. Then there are Seifert bundles 
$f:L\to ((k+1)\#\c\p^2,\Delta)$ 
such that $L$ is simply connected, $i(L)=1$ and
$$
H_2(L,\z)=\z^k+\sum_{p,i} \bigl(\z/p^i\bigr)^{c(p^i)}.\qed
$$
\end{cor}

\begin{rem} Note that a surface $D_{ij}$ with genus 0
does not contribute to $H_2(L,\z)$. Thus if we add
new $m_{ij}$ with $g_{ij}=0$ to our collection, we get
the same total space for the Seifert bundle.
The number and genera of the $D_{ij}$ are easy to determine
from the subset of $M$ where the action is not free.
Thus for each such $M$ we get infinitely many
 topologically distinct circle actions.
\end{rem}

\section{Quasi--regular contact and Sasakian structures}

 An  interesting case
of $S^1$-actions arises when $(M,\eta)$ is a contact manifold
and the Reeb vector field gives an $S^1$-action. These are called
{\it quasi--regular} contact structures. These are distinguished
by the fact that $d\eta$ descends to an (orbifold) symplectic structure
on $X=M/S^1$ and the $D_i\subset X$ are symplectic suborbifolds.

I can not say anything useful about the general
contact case, but as a first step one may consider
Sasakian structures (see \cite{bg-review} for a recent survey).
 For our purposes, these are 
 Seifert bundles over algebraic orbifolds.
That is, Seifert bundles $f:M\to (X,\sum (1-\frac1{m_i})D_i)$
where 
 $X$ is a complex  algebraic (possibly singular) surface and 
$D_i\subset X$ are
 complex  algebraic curves.

Since symplectic 4--manifolds are close to
complex algebraic surfaces in many respect,
one may hope that some features of the Sasakian case
continue to hold for contact manifolds as well.

 \cite[Cor.81]{ko-es5} shows that not every
simply connected 
rational homology sphere
admits a Sasakian structure, but now we see
that the restrictions found there
 are purely topological. That is, they are obstructions to the
existence of  a fixed point free circle action as well.

Here is a relatively simple result which shows that
Sasakian structures impose additional topological restrictions
beyond those that come from the circle  action itself
(\ref{main.thm}).

\begin{lem}\label{sasak.rests}
 Let $M\to (X,\sum (1-\frac1{m_i})D_i)$ be a Seifert bundle
over an algebraic orbifold.
Assume that $H_1(M,\z)=0$ and 
$$
H_2(M,\z)=\sum_{p,i} \bigl(\z/p^i\bigr)^{c(p^i)}.
$$
Then there is a degree 2 polynomial $q$ with integer coefficients
such that $q(\z)$  contains all but 10 elements of 
 the set 
$\{c(p^i)\}$.

In particular,  $\{c(p^i)\}$ 
 contains at most $12+2\sqrt{N}$ elements 
of any interval of length $N$.
\end{lem}

\begin{exmp} 
Let the $p_i$ be different prime numbers and 
$M\to (X,\sum (1-\frac1{m_i})D_i)$  a Seifert bundle
over an algebraic orbifold such that $H_1(M,\z)=0$ and 
$$
H_2(M,\z)=\sum_{i=1}^k \bigl(\z/p_i\bigr)^{2i}.
$$
Then (\ref{sasak.rests}) implies  that $k\leq  23$.

On the other hand, 
the conditions of (\ref{main.thm}) are satisfied for any $k$,
and so there are such Seifert bundles 
for any $k$ and $p_i$ even over $X=\c\p^2$.
The surfaces $D_i$, however,  can not be chosen complex
algebraic for $k>23$.

These give  examples  of simply connected
5--manifolds which admit a fixed point free circle action
yet have no Sasakian structure.

With more careful estimates and some case analysis,
 one should be able to reduce 23 to about 10, but it
may be hard to get a sharp result.
\end{exmp}

Proof. If $X$ is an algebraic surface with quotient singularities,
the Bogomolov--Miyaoka--Yau--type inequalities of
\cite[10.8, 10.14]{k-etal} and \cite[9.2]{ke-mc} imply that 
$$
{\textstyle \sum_x}\bigl(1-\tfrac1{r_x}\bigr)<e(X),
$$
where $r_x$ is the order of the fundamental group of the link of 
 a singular point $x\in X$
and $e(X)$ is the topological Euler number.
In our case $e(X)=3$ and so 
$X$ has at most 5 singular points.
Thus by (\ref{dim5.nice.say})
 there are at most 10 curves $D_i$ which pass through
 a singular point.

Once $D_i$ is contained in the smooth locus of $X$,
its genus is computed by the adjunction formula
$$
2g(D_i)=(D_i\cdot (D_i+K_X))+2.
$$
Since $X$ has Picard number 1,  we conclude that
there is a  degree 2 polynomial with integer coefficients $q(t)$
such that $2g_i(D)\in q(\z)$.

Thus we are finished by the following easy lemma:

\begin{lem} Let $q(t)=at^2+bt+c$ with $a>0$.
 Then
the set $q(\z)$ intersects every interval of length $N$
in at most $2+2\sqrt{N/a}$ elements.\qed
\end{lem}

\begin{exmp} Let $X=\p^2$ and $D_i\subset X$ be a smooth curve
of degree $i$ for $i=3,\dots, n$, hence
$g(D_i)=\binom{i-1}{2}$. Choose pairwise relatively prime
integers $m_i$.  By (\ref{seif.exist.thm}) 
and (\ref{H_2.tors.onto})    there is
a simply connected Seifert bundle
$M\to (\c\p^2, \sum (1-\frac1{m_i})D_i)$ such that  
$$
H_2(M,\z)=\sum_{i=3}^n \bigl(\z/m_i\bigr)^{2\binom{i-1}{2}}.
$$
This gives  $\sqrt{N}-1$ different values $c_{ip}$
in the interval $[1,N]$ for $N=n^2$.
\end{exmp}

\begin{ack} I thank A.\ Adem, Ch.\ Boyer, K.\ Galicki and  Z.\ Szab\'o
for  useful conversations and e-mails.
Research
was partially supported by the NSF under grant number
DMS-0200883. 
\end{ack}

\bibliography{refs}

\vskip1cm

\noindent Princeton University, Princeton NJ 08544-1000

\begin{verbatim}kollar@math.princeton.edu\end{verbatim}

\end{document}